\documentclass[aos,MSNbibl,seceqn,citesort,dvips]{arximspdf}
\usepackage{graphicx}

\doi{10.1214/11-AOS951}
\volume{39}
\issue{6}
\pubyear{2011}
\firstpage{3417}
\lastpage{3440}

\makeatletter

\newtheorem{theorem}{Theorem}[section]
\newtheorem{proposition}[theorem]{Proposition}
\newtheorem{lemma}[theorem]{Lemma}

\newproclaim{remark}[theorem]{Remark}

\renewcommand{\O}{{\mathrm{O}}}
\renewcommand{\P}{{\mathcal{P}_2}}
\renewcommand{\L}{\mathcal{L}}
\newcommand{\M}{\mathcal{M}}

\newcommand{\tr}{\operatorname{tr}}
\newcommand{\GL}{\operatorname{GL}}
\newcommand{\mE}{\mathbb{E}}
\newcommand{\R}{\mathbb{R}}
\newcommand{\C}{\mathbb{C}}
\newcommand{\F}{{{\mathcal F}_\varphi}}
\newcommand{\He}{{\mathcal H}}
\newcommand{\I}{{\mathbf I}}
\newcommand{\Pm}{{\mathcal{P}_m}}
\newcommand{\Pos}{{\mathcal P}}
\newcommand{\fhat}{\hat{f}}
\newcommand{\what}{\hat{w}}

\newcommand{\im}{\operatorname{Im}}
\newcommand{\re}{\operatorname{Re}}
\newcommand{\half}{{\frac1 2}}
\newcommand{\D}{{{\mathcal D}}}

\newcommand{\argmin}{\mathop{\arg\min}}
\newcommand{\diag}{\operatorname{diag}}
\newcommand{\MISE}{\operatorname{MISE}}

\makeatother

\begin{document}
\begin{frontmatter}

\title{Minimax estimation for mixtures of~Wishart~distributions}
\runtitle{Minimax estimation for Wishart mixtures}

\begin{aug}
\author[A]{\fnms{L. R.} \snm{Haff}\ead[label=e1]{haff@ucsd.edu}},
\author[B]{\fnms{P. T.} \snm{Kim}\corref{}\thanksref{t2}\ead[label=e2]{pkim@uoguelph.ca}},
\author[C]{\fnms{J.-Y.} \snm{Koo}\thanksref{t3}\ead[label=e3]{jykoo@korea.ac.kr}}
\and
\author[D]{\fnms{D. St. P.} \snm{Richards}\thanksref{t4}\ead[label=e4]{richards@stat.psu.edu}}
\runauthor{Haff, Kim, Koo and Richards}
\affiliation{University of California at San Diego,
University of Guelph, Korea University and Penn State University}
\address[A]{L. R. Haff\\
Department of Mathematics \\
University of California at San Diego \\
La Jolla, California 92093\\
USA \\
\printead{e1}}
\address[B]{P. T. Kim\\
Department of Mathematics\hspace*{35pt}\\
\quad and Statistics \\
University of Guelph \\
Guelph, Ontario N1G 2W1\\
Canada \\
\printead{e2}}
\address[C]{J.-Y. Koo\\
Department of Statistics \\
Korea University \\
Anam-Dong Sungbuk-Ku \\
Seoul, 136-701\\
Korea \\
\printead{e3}}
\address[D]{D. St. P. Richards\\
Department of Statistics\\
Penn State University \\
University Park, Pennsylvania 16802\\
USA \\
\printead{e4}} 
\end{aug}

\thankstext{t2}{Supported in part by the Natural Sciences
and Engineering Research Council of Canada, Grant DG 46204.}

\thankstext{t3}{Supported in part by Basic Science Research Program
through the
National Research Foundation of Korea (NRF) funded by the Ministry of
Education, Science and Technology (2011-0027601).}

\thankstext{t4}{Supported in part by NSF
Grant DMS-07-05210.}

\received{\smonth{8} \syear{2009}}
\revised{\smonth{7} \syear{2011}}

%
\begin{abstract}
The space of positive definite symmetric matrices has been studied
extensively as a means of understanding dependence in multivariate data
along with the accompanying problems in statistical inference. Many
books and papers have been written on this subject, and more recently
there has been considerable interest in high-dimensional random
matrices with particular emphasis on the distribution of certain
eigenvalues. With the availability of modern data acquisition
capabilities, smoothing or nonparametric techniques are required that
go beyond those applicable only to data arising in Euclidean spaces.
Accordingly, we present a Fourier method of minimax Wishart mixture
density estimation on the space of positive definite symmetric
matrices.
\end{abstract}

%

\end{frontmatter}

\section{Introduction}
\label{introduction}

The space of positive definite symmetric matrices has been studied
extensively in statistics as a means of understanding dependence in
multivariate data along with the accompanying problems in statistical
inference. Many books and papers, for example,
\cite{haff1,haff2,haff3,loh,Muirhead,stein-rietz} and \cite{stein},
have been written on this subject, and there has been considerable
interest recently in high-dimensional random matrices with particular
emphasis on the distribution of certain eigenvalues \cite{johnstone08}
and on graphical models \cite{lm}.\looseness=1

In this paper we consider the problem of estimating the mixing
density of a~continuous mixture of Wishart distributions.
We construct a nonparametric estimator of that density and
obtain minimax rates of convergence for the estimator.
Throughout this work, we adopt, as a guide,
results developed for the classical problem of deconvolution density
estimation on Euclidean spaces; see, for example,
\cite{diggle,fan,Koo,mr,butsy} and \cite{Zhang}.
Much of the difficulty with the space of positive definite symmetric
matrices is due to the fact that mathematical analysis on the space
is technically demanding.
Helgason \cite{he2} and Terras \cite{Terras-II} provide much insight
and technical innovation,
however, and we make extensive use of these methods.

We summarize the paper as follows. In Section \ref{secwishart} we
discuss and set up the notation for Wishart mixtures. In Section
\ref{secdecon} we begin by reviewing the necessary Fourier methods
which allow us to construct a nonparametric estimator of the mixing
density, and then we provide the estimator. The minimax property of our
nonparametric estimator is stated in Section \ref{secmain} along with
supporting results. Section~\ref{secNS} presents simulation results as
well as an application to finance examining real financial data.
Finally, Sections~\ref{secproof-upper} and \ref{secproofs-lower}
present the proofs.

\section{Wishart mixtures}
\label{secwishart}

Throughout the paper, for any square matrix $y$, we denote the
trace and determinant of $y$ by $\tr y$ and $|y|$, respectively;
further, we denote by $\I_m$ the $m \times m$ identity matrix.
We will denote by $\Pm$ the space of $m \times m$ positive
definite symmetric matrices.

For $s = (s_1,\ldots,s_m) \in\C^m$ with $\re(s_j+\cdots+s_m) >
(j-1)/2$, $j=1,\ldots,m$, the multivariate gamma function is
defined as
%
%
\begin{equation}
\label{prod-gamma}
\Gamma_m(s_1,\ldots,s_m) = \pi^{m(m-1)/4}
\prod_{j=1}^m \Gamma\biggl(s_j+\cdots+s_m-\frac12(j-1)\biggr),
\end{equation}
where $\Gamma(\cdot)$ denotes the classical gamma function.

We denote by $G$ the general linear group $\GL(m,\R)$ of all
$m \times m$ nonsingular real matrices, by $K$ the group
$\O(m)$ of $m \times m$ orthogonal matrices and
by $A$ the group of diagonal positive definite matrices.
The group $G$ acts transitively on $\Pm$ by the action
%
%
\begin{equation}
\label{gp-action}
G\times\Pm\to\Pm,\qquad (g,y) \mapsto g'yg,
\end{equation}
$g \in G$, $y \in\Pm$, where $g'$ denotes the transpose of $g$.
Under this group action, the isotropy group of the identity in $G$
is $K$; hence the homogeneous space~$K \setminus G$ can be
identified with $\Pm$ by the natural mapping from
$K\setminus G\to\Pm$ that sends $Kg \mapsto g'g$.
In distinguishing between left and right cosets, we place the
quotient operation on the left and right of the group, respectively.

For $y = (y_{ij}) \in\Pm$, define the measure
\[
{\mathrm{d}_*}y = |y|^{-(m+1)/2} {\prod_{1 \le i \le j \le m} d y_{ij}}.
\]
It is well known that the measure ${\mathrm{d}_*}y$ is invariant under
action (\ref{gp-action}). Relative to the dominating measure
${\mathrm{d}_*}y$, the probability density function of the standard Wishart
distribution with $N$ degrees of freedom is
%
%
\begin{equation}
\label{wishart}
w(y) = \frac{1}{2^{Nm/2}\Gamma_m(0,\ldots,0,N/2) }
|y|^{N/2} \exp\biggl(-\frac12 \tr y\biggr),
\end{equation}
$y \in\Pm$. Consequently, for $\sigma\in\Pm$, we note
that $\tr(\sigma^{-1/2} y \sigma^{-1/2}) = \tr(\sigma^{-1} y)$
and $|\sigma^{-1/2} y \sigma^{-1/2}| = |\sigma^{-1} y|$. It then
follows that, relative to the dominating measure~${\mathrm{d}_*}y$, the
density of the general Wishart distribution, with covariance
parameter~$\sigma$, is $w(\sigma^{-1}y)$, $y \in\Pm$.

Suppose next that $\sigma$ is a random matrix
and, relative to the dominating measure~${\mathrm{d}_*}\sigma$, has a
continuous mixing density, $f$, that is invariant under the
action~(\ref{gp-action}). By integration with respect to
$\sigma$, the continuous Wishart mixture density is given by
%
%
\begin{equation}
\label{marginal}
r(y) = \int_{\Pm} f(\sigma) w (\sigma^{-1}y )
\,{\mathrm{d}_*}
\sigma,
\end{equation}
$y \in\Pm$.
For the case in which $m=1$, the standard Wishart density is
essentially a chi-square density, in which case (\ref{marginal})
is a continuous mixture of chi-square densities.

In general, (\ref{marginal}) is a convolution operation for functions
on $\Pm$.
We denote by $x^{1/2}$ any matrix with $x^{t/2} x^{1/2}=x$,
where $x^{t/2} = (x^{1/2})'$ and denote $x^{-t/2} = (x^{-1/2})^t$.
Define $X \circ Z$, the convolution of two random matrices $X$ and $Z$
which are distributed on $\Pm$ by
\[
X \circ Z = X^{t/2} Z X^{1/2}
\]
and $f_1 * f_2$, the convolution of $f_1 \in L^1(\Pm)$ and $f_2 \in
L^1(\Pm)$ by
\[
(f_1 * f_2)(y) = \int_\Pm f_1(x) f_2 (x^{-t/2} y x^{-1/2} )
\,{\mathrm{d}_*}x
\qquad\mbox{for } y \in\Pm,
\]
where $L^q(\Pm)$ is the space of integrable functions raised to the
$q$th power on $\Pm$
for $q \geq1$.
If $X$ and $Z$ with densities $f$ and $w$, respectively, are independent,
then $Y = X \circ Z$ has the density $r = f * w$
since $w(\sigma^{-t/2} y \sigma^{-1/2}) = w(\sigma^{-1}y)$.
Finally, (\ref{marginal}) can be transformed into a scalar multiplication;
see Section \ref{ssechfc}.

\section{Fourier analysis on $\Pm$ and
estimation of the mixing density}
\label{secdecon}

In this section we review the Fourier methods needed to transform
the convolution product (\ref{marginal}) and to construct a
nonparametric estimator of the mixing density $f$.

\subsection{The Helgason--Fourier transform}
\label{sechf}

For $y \in\Pm$, denote by $|y_j|$
the principal minor of order $j$, $j=1,\ldots,m$.
For
$s \in\C^m$, the \textit{power function}~$p_s \dvtx\allowbreak \Pm\to\C$ is
%
%
\begin{equation}
\label{powerfunct}
p_s(y) = \prod_{j=1}^m |y_j|^{s_j},
\end{equation}
$y \in\Pm$.
Let ${\mathrm{d}_*}k$ denote the Haar measure on $K$, normalized to have total
volume equal to one; then
%
%
\begin{equation}
\label{zonalsphfunct}
h_s(y) = \int_K p_s(k' yk) \,{\mathrm{d}_*}k,
\end{equation}
$y \in\Pm$, is the \textit{zonal spherical function} on $\Pm$.
It is well known that the functions $h_s$ are fundamental to
harmonic analysis on symmetric spaces \cite{he2,Terras-II}.
If $s_1,\ldots,s_m$ are nonnegative integers
then, up to a constant factor, (\ref{zonalsphfunct}) is an
integral formula for the zonal polynomials which arise in many
aspects of multivariate statistical analysis
\cite{Muirhead}, pages 231 and 232.

Let $C_c^{\infty}(\Pm)$ denote the space of infinitely
differentiable, compactly supported, complex-valued functions
$f$ on $\Pm$; also, let
\[
C_c^\infty(\Pm/K) = \{f \in C_c^\infty(\Pm)\dvtx
f(k'yk)=f(y) \mbox{ for all } k \in K, y \in\Pm\}.
\]
For $s \in\C^m$ and $k \in K$, the
\textit{Helgason--Fourier transform} (\cite{Terras-II}, page 87) of
a~function $f \in C_c^\infty(\Pm)$ is
%
%
\begin{equation}
\label{hf-transform}
\He f(s,k) = \int_{\Pm} f(y)
{\overline{p_s(k'yk)}} \,{\mathrm{d}_*}y,
\end{equation}
where ${\overline{p_s(k'yk)}}$ denotes complex conjugation.

For the case in which $f \in C_c^\infty(\Pm/K)$, we make the
change of variables $y \mapsto k_1'yk_1$ in (\ref{hf-transform}),
$k_1 \in K$, and integrate with respect to the Haar measure
${\mathrm{d}_*}k_1$. Applying the invariance of $f$ and formula
(\ref{zonalsphfunct}), we deduce that~$\He f(s,k)$ does not
depend on $k$. Specifically, $\He f(s,k) = \fhat(s)$ where
%
%
\begin{equation}
\label{hfs-transform}
\fhat(s) = \int_{\Pm} f(y)
{\overline{h_s(y)}} \,{\mathrm{d}_*}y,
\end{equation}
$s \in\C^m$, is the \textit{zonal spherical transform} of $f$.

In the case of the standard Wishart density (\ref{wishart}),
which is a $K$-invariant function, the zonal spherical
transform is well known (Muirhead \cite{Muirhead}, pa\-ge~248;
Terras~\cite{Terras-II}, pages 85 and 86):
\[
\what(s) = \frac{\Gamma_m(s_{m-1},\ldots,s_1,-(s_1+\cdots+s_m)+N/2)}
{\Gamma_m(0,\ldots,0,N/2)} h_s\biggl(\frac12 \I_m\biggr).
\]

\subsection{The convolution property of the Helgason--Fourier transform}
\label{ssechfc}

The following result shows that the convolution operation can be
transformed into a scalar multiplication.\vadjust{\goodbreak}
%
%
\begin{proposition} \label{propHFconvC}
Suppose $X$ and $Z$ with densities $f_X$ and $f_Z$, respectively, are
independent, and $Z$ is $K$-invariant. Let $f_Y$ be the density of $Y =
X \circ Z$. Then
\[
\He f_Y(s,k) = \He f_X(s,k) \hat f_Z(s)
\qquad\mbox{for } s \in\C^m \mbox{ and } k \in K.
\]
\end{proposition}
\begin{pf}
Note
\[
\He f_Y(s,k)
= \mE p_{\bar s}(k' Y k)
= \mE p_{\bar s}(k' X^{t/2} Z X^{1/2} k).
\]
Using the $KAN$-Iwasawa decomposition of $X^{t/2} k$ (Terras
\cite{Terras-II}, page 20),
we have $X^{1/2} k = H U$ for $H \in K$ and $U$, an upper triangular matrix.
Observe
\begin{eqnarray*}
\mE p_{\bar s}(k' X^{t/2} Z X^{1/2} k)
&=& \mE_X \{ p_{\bar s}(U' U)
\mE_{Z|X} p_{\bar s}(H'Z H) \} \\
&=& \hat f_Z(s) \mE_X \{ p_{\bar s} (U' U)
\} \\
&=& \He f_X(s,k) \hat f_Z(s),
\end{eqnarray*}
where Proposition 1 of Terras \cite{Terras-II}, page 39, is used for
the first equality.
\end{pf}

\subsection{The inversion formula for the Helgason--Fourier transform}
\label{ssechfi}

For $a_1,\allowbreak a_2 \in\C$ with $\re(a_1), \re(a_2) > 0$, let
\[
B(a_1,a_2) = \frac{\Gamma(a_1) \Gamma(a_2)}{\Gamma(a_1 + a_2)}
\]
denote the classical beta function. For
$s \in\C^m$ such that
$\re(s_i+\cdots+s_j) > -\frac12(j-i+1)$ for all
$1 \le i<j \le m-1$, the \textit{Harish--Chandra} $c$-\textit{function} is
%
%
\begin{equation}\label{cfunction}
c_m(s) = \prod_{1 \le i < j \le m-1}
\frac{B(1/2,s_i+\cdots+s_j+(j-i+1)/2)}
{B(1/2,(j-i+1)/2)} .
\end{equation}
Let $\rho\equiv(\frac12,\ldots,\frac12,\frac
14(1-m))$,
and set
%
%
\begin{eqnarray}
\label{omega-m}
\omega_m &=& \frac{\prod_{j=1}^m \Gamma(j/2)}
{(2\pi i)^{m} \pi^{m(m+1)/4} m!} ,
\\
%
%
\label{Cmrho}
\C^m(\rho) &=& \{s \in\C^m\dvtx \re(s) = -\rho\}
\end{eqnarray}
and
\[
{\mathrm{d}_*}s = \omega_m |c_m(s)|^{-2} \,{\mathrm{d}_*}s_1 \cdots{
\mathrm{d}_*}s_m.
\]

Let $M = \{{\operatorname{diag}}(\pm1,\ldots,\pm1)\}$ be the set of
$m \times m$ diagonal matrices with entries $\pm1$ on the diagonal;
then $M$ is a subgroup of $K$ and is of order $2^m$.
By factorizing the Haar measure ${\mathrm{d}_*}k$ on $K$, it may be shown
(\cite{Terras-II}, page 88) that there exists an invariant measure
${\mathrm{d}_*}\bar{k}$ on the coset space $K/M$ such that
\[
\int_{\bar{k} \in K/M} {\mathrm{d}_*}\bar{k} = 1.\vadjust{\goodbreak}
\]

The \textit{inversion formula for the Helgason--Fourier transform}
$\He$ in (\ref{hfs-transform}) is that if
$f \in C^\infty_c(\Pm)$, then \cite{he2,Terras-II}
%
%
\begin{equation}
\label{hf-inversion}
f(y) = \int_{\C^m(\rho)}
\int_{\bar{k} \in K/M} \He f(s,k) p_s(k'yk) \,{\mathrm{d}_*}\bar{k}\,
{\mathrm{d}_*}s,
\end{equation}
$y \in\Pm$. In particular, if $f \in C^\infty_c(\Pm/K)$, then
\[
f(y) = \int_{\C^m(\rho)} \fhat(s) h_s(y) \, {\mathrm{d}_*}s,
\]
$y \in\Pm$, and there also holds the \textit{Plancherel formula},
%
%
\begin{equation}
\label{hf-plancherel}
\int_{\Pm} |f(y)|^2 \,{\mathrm{d}_*}y = \int_{\C^m(\rho)} \int
_{K/M}
|\He f(s, k)|^2 \,{\mathrm{d}_*}\bar k \,{\mathrm{d}_*}s .
\end{equation}
We refer to Terras \cite{Terras-II}, page 87 ff., for full details
of the inversion formula and for references to the literature.

\subsection{Eigenvalues, the Laplacian and Sobolev spaces}
\label{ssecels}

For $y = (y_{ij}) \in\Pm$, we define the $m \times m$ matrix
of partial derivatives,
\[
\frac{\partial}{\partial y} = \biggl(\frac12(1+\delta
_{ij})\,
\frac{\partial}{\partial y_{ij}}\biggr),
\]
where $\delta_{ij}$ denotes Kronecker's delta. The Laplacian,
$\Delta$, on $\Pm$ can be written (\cite{Terras-II}, page 106)
in terms of the local coordinates $y_{ij}$ as
\[
\Delta= -\tr\biggl(\biggl(y\,\frac{\partial}{\partial y}
\biggr)^2\biggr).
\]

The power function $p_s$ in (\ref{powerfunct}) is an
eigenfunction of $\Delta$ (see \cite{Muirhead}, page 229,
\cite{Richards1}, page 283, \cite{Terras-II}, page 49). Indeed,
let $r_j = s_j+s_{j+1}+\cdots+s_m+\frac14(m-2j+1)$,
$j=1,\ldots,m$, and define
%
%
\begin{equation}
\label{eigenvalue}
\lambda_s = -(r_1^2 + \cdots+ r_m^2 ) + \tfrac{1}{48}m(m^2
- 1);
\end{equation}
then $\Delta p_s(Y) = \lambda_s p_s(Y)$. Since $\re(s) = -\rho$
then each $r_j$, $j=1,\ldots,m$, is purely imaginary; hence,
$\lambda_s > 0$, $s \in C_m(\rho)$.

The operator $\He$ changes the effect of invariant differential
operators on functions to pointwise multiplication: if
$f \in C^\infty_c(\Pm)$, then
\[
\He(\Delta f)(s,k) = \lambda_s \He f(s,k),
\]
$s \in\C^m$, $k \in K$ (\cite{Terras-II}, page 88). For
$\varphi> 0$, we therefore define the fractional power,
$\Delta^{\varphi/2}$, of $\Delta$, as the operator such that
\[
\He(\Delta^{\varphi/2} f)(s,k) = \lambda_s^{\varphi/2} \He f(s,k),
\]
$f \in C^\infty_c(\Pm)$. Having constructed $\Delta^{\varphi/2}$,
we define the \textit{Sobolev class},
\[
\F= \{f \in C^{\infty}(\Pm) \dvtx
\|\Delta^{\varphi/2}f\|^2 < \infty\},\vadjust{\goodbreak}
\]
where for $f \in C^\infty(\Pm)$,
\[
\|f\| = \biggl(\int_{\Pm} |f(y)|^2 \,{\mathrm{d}_*}y \biggr)^{1/2}
\]
denotes the $L^2(\Pm)$-norm with respect to the measure
${\mathrm{d}_*}y$. For $Q > 0$, we also define the \textit{bounded Sobolev class},
\[
\F(Q) = \{f \in C^{\infty}(\Pm) \dvtx
\|\Delta^{\varphi/2}f\|^2 < Q \}.
\]

\section{Main result}
\label{secmain}

In this section we will present the main result.
We do so by applying the Helgason--Fourier transform to the mixture
density (\ref{marginal})
so that
%
%
\begin{equation}\label{hc-convoprop-prior}
\He r(s,k) = \He f (s,k) \what(s),
\end{equation}
$s \in\C^m$, $k \in K$; see Proposition \ref{propHFconvC}.
Having observed a random sample~$Y_1,\allowbreak\ldots,Y_n$ from the mixture
density, $r$, in
(\ref{marginal}), we estimate $\He r (s,k)$ by its
\textit{empirical Helgason--Fourier transform},
%
%
\begin{equation}
\label{emp}
{\He_n r}(s,k)= \frac{1}{n} \sum_{\ell=1}^n
\overline{p_s(k'Y_{\ell}k)}.
\end{equation}
On substituting (\ref{emp}) in (\ref{hc-convoprop-prior}), together
with the assumption that $\what(s) \neq0$, $s \in\C^m$, we obtain
\[
{\He_n f }(s,k) = \frac{\He_n r(s,k)}{{\what}(s)},
\]
$s \in\C^m$, $k \in K$.

Analogous with classical Euclidean deconvolution, we introduce a
smoothing parameter $T = T(n)$ where $T(n) \to\infty$ as
$n \to\infty$, and then we apply the inversion formula
(\ref{hf-inversion}) using a spectral cut-off based on the
eigenvalues of~$\Delta$. First, we introduce the notation
\[
\C^m(\rho,T) = \{s \in\C^m(\rho)\dvtx \lambda_s < T \},
\]
where $\C^m(\rho)$ is defined in (\ref{Cmrho}). We now define
%
%
\begin{equation}\label{emp-hf-inversion}
f_n(y) = \int_{\C^m(\rho,T)} \int_{\bar{k} \in K/M}
\frac{\He_n r(s,\bar{k})}{{\what}(s)}
p_s(\bar{k}'y\bar{k}) \,{\mathrm{d}_*}\bar{k} \,{\mathrm{d}_*}s,
\end{equation}
$y \in\Pm$, and take this as our
nonparametric estimator of $f$.

We now state the minimax result for the estimator (\ref{emp-hf-inversion}).
Let $C$ denote a~generic positive constant independent of $n$.
For two sequences of real numbers
$\{a_n\}$ and~$\{b_n\}$, we use the
notations $a_n \ll b_n$ and $a_n \gg b_n$
to mean $a_n < C b_n$ and $a_n > C b_n$, respectively,
as $n \to\infty$.
Moreover, $a_n \asymp b_n$ means that $a_n \ll b_n$ and $a_n \gg b_n$.
%
%
\begin{theorem}
\label{thmwishart1}
Suppose $f$ is a density on $\Pm$ and $N > (m-1)/2$.
Then, for the Wishart mixture (\ref{marginal}),
%
%
\begin{equation}\label{mainresult}
\sup_{f \in\F(Q)}
\mE\|f_n - f\|^2 \ll(\log n)^{-2\varphi}
\end{equation}
and for any estimator $g_n$ of $f$,
%
%
\begin{equation}\label{eqkower}
\inf_{g_n} \sup_{f \in\F(Q)}
\mE\|g_n - f\|^2 \gg(\log n)^{-2\varphi}.
\end{equation}
\end{theorem}

We now provide some comments about this result.
In the situation whe\-re~(\ref{marginal}) is a finite sum, so that
\[
r(y) = \sum_{\ell=1}^q f_\ell w(\sigma_\ell^{-1}y)
\quad\mbox{and}\quad \sum_{\ell=1}^q f_\ell=1,
\]
we have the finite mixture model. Methods for recovering the
mixing coefficients can be covered by the techniques employed in
\cite{chentan}. We note that the continuous mixture
model is a generalization of this approach.

It is noted that the condition $f$ is a density and seems to be mild.
The upper bound property of (\ref{mainresult}) is established in
\cite{kr2}, Theorem 3.3, with $\beta=1/2$.
In the latter, the moment condition
%
%
\begin{equation}
\label{moment}
\int_{\Pm} |y_1|^{-1} \cdots|y_{m-1}|^{-1} |y|^{(m-1)/2} r(y)
\,{\mathrm{d}_*}
y < \infty,
\end{equation}
on the principal minors $|y_1|,\ldots,|y_m|$ of $y \in\Pm$ is
assumed. In our theorem, we did not impose this moment condition as
condition (\ref{moment}) is automatically satisfied. This is pointed out
and commented upon below in the proof.

To derive the lower bound for estimating $f$ in the
$L^2(\Pm)$-norm, we shall follow the standard Euclidean
approach. Thus we choose a pair of functions~$f^0$, $f^n \in\F(Q)$, and, with $w$ denoting the Wishart
density (\ref{wishart}), we shall show that, for some constants
$C_1, C_2 > 0$,
\[
\|f^n - f^0 \|^2 \geq C_1 (\log n)^{-2\varphi}
\]
and
%
%
\begin{equation}\label{chisquaredist}
\chi^2(f^0*w,f^n*w) \leq\frac{C_2} n,
\end{equation}
where
\[
\chi^2(g_1, g_2) = \int_{\Pm} \frac{(g_1(y) - g_2(y)
)^2}{g_1(y)} \,{\mathrm{d}_*}y.
\]

Precisely, let us suppose we can choose $f^0 \in\F(Q)$
and a perturbation $\psi\in\F(Q)$, and, for
$\delta= \delta_n > 0$, let $\psi^{\delta}$
be a scaling of $\psi$ such that
$\|\psi^{\delta}\| \asymp\delta^{-1/2}\|\psi\|$.
Define
\[
f^n = f^0 + C_\psi\delta^{-\varphi+ 1/2} \psi^{\delta}.
\]
If $\delta$ can be chosen so that
\[
\chi^2(f^0*w,f^n*w) \le Cn^{-1},
\]
then the lower bound rate of convergence is
determined by $\delta^{-2\varphi}$. We shall develop such a
construction and, moreover, do so in a way such that
$\delta\asymp\log n$ as $n \to\infty$.
%
%
\begin{remark}
The profound influence of Charles Stein on covariance estimation
originates largely from his Rietz lecture; see \cite{stein-rietz}. The
idea is that for certain loss functions over $\Pm$, the usual estimator
of the covariance matrix parameter is inadmissible. Through an unbiased
estimation of the risk function over covariance matrices, Stein was
able to improve upon the usual estimator by pooling the observed
eigenvalues of the sample covariance matrix. Subsequent to this,
through a series of papers, improvements were obtained in Haff
\cite{haff1,haff2,haff3}. Other related works include Takemura
\cite{takemura}, Lin and Perlman \cite{linperl} and Loh \cite{loh}, to
name a few.

In this paper, we contribute to the case in which one observes
data from a continuous mixture of Wishart distributions, not merely
a sample from a~single distribution. Therefore, the parameter of interest
would be the mixing density of the covariance parameters.
And the nonparametric estimator of the mixing density (\ref{emp-hf-inversion})
is an attractive candidate because of its minimax property.
Based on this procedure, one could consider
the moment, or mode, of $f_n$, as a possible estimator of the
corresponding population parameters. Alternatively, one could take a
nonparametric empirical Bayes approach as in Pensky \cite{pensky}.
\end{remark}

\section{Numerical aspects and an application to finance}
\label{secNS}

This section presents numerical aspects for the $m = 2$ case with an
application to finance.

\subsection{Computation of estimators}
\label{subsecCS}

Suppose $X$ and $Z$ are independent\break with~$Z$ having a Wishart distribution.
Let $Y = X^{t/2} Z X^{1/2}$ where $X^{1/2}$ is upper triangular.
For visualization, we display estimators of the marginal density for~$D$
where $X = H' D H$ with $H \in K$ and $D \in A$.
Let
\[
\hat r_n(s) = \frac1 n \sum_{j=1}^n \overline{h_s(E_j)}
\quad\mbox{and}\quad
\hat f_n(s) = \frac{\hat r_n(s)} {\hat w(s)},
\]
where $E_j \in A_+$ denotes the diagonal matrix
of eigenvalues of $Y_j$, $j = 1, \ldots, n$.
Denote by $f^D$ the density of eigenvalues of $X$.
Then, the estimator for $f^D$ is given by
%
%
\begin{equation}\label{eqfXhat2}\qquad
f^D_n(a) = \int_{\C^2(\rho, T)} \re\{\hat f_n(s) h_s(a)
\}\,
{\mathrm{d}_*}s\qquad
\mbox{for } a = \diag(a_1, a_2) \in A.
\end{equation}

Consider the computation of $h_s(a)$
when
\[
s = -\rho+ ib = (-1/2+ib_1, 1/4+ib_2)\vadjust{\goodbreak}
\]
so that $\re(s) = -\rho$.
From pages 90 and 91 of \cite{Terras-II}, the spherical function is
given by
\[
h_s(a) = (a_1 a_2)^{i(b_2+b_1/2)}
P_{-1/2+ i b_1}\bigl(\cosh\bigl(\log\bigl(\sqrt{a_1/a_2}
\bigr)\bigr)\bigr),
\]
where Legendre function $P_{-1/2+ i t}(x)$
can be computed using \texttt{conicalP\_0(t, x)} in the \texttt{gsl}
library in \texttt{R}.

The mean integrated squared error (MISE) of $f^D_n$ is defined by
\[
\MISE(f_n) = \mE\int_A \bigl(f^D_n(a) - f^D(a) \bigr)^2 \,{\mathrm{d}_*}a.
\]
It is reasonable to choose $T$ which minimizes $\MISE(f^D_n)$
or equivalently
\[
M(T) = \mE\int_A (f^D_n(a))^2\, {\mathrm{d}_*}a - 2 \mE\int
_A f^D_n(a)
f^D(x)\, {\mathrm{d}_*}a.
\]
One can find an unbiased estimator $M_0(T)$ of $M(T)$; see \cite{kkkoz}.
We choose $\hat T$ by
\[
\argmin_T M_0(T).
\]
Monte Carlo approximation is used for integration of (\ref{eqfXhat2})
and $M_0(T)$.

\subsection{Simulation}
\label{subsecsimul}

Denote by $W_N(\sigma)$ the Wishart density with degrees of freedom $N$
and covariance matrix $\sigma$.
We generate data as follows.
For $j = 1, \ldots, n$:
\begin{itemize}
\item
generate $Z_j \sim W_{20}(I_2)$;
\item
generate $X_j \sim f$;
\item
do a Cholesky decomposition of $X_j = (X_j)^{t/2} (X_j)^{1/2}$, and
calculate $ Y_j = (X_j)^{t/2} Z_j (X_j)^{1/2}$.
\end{itemize}

%
\begin{figure}

\includegraphics{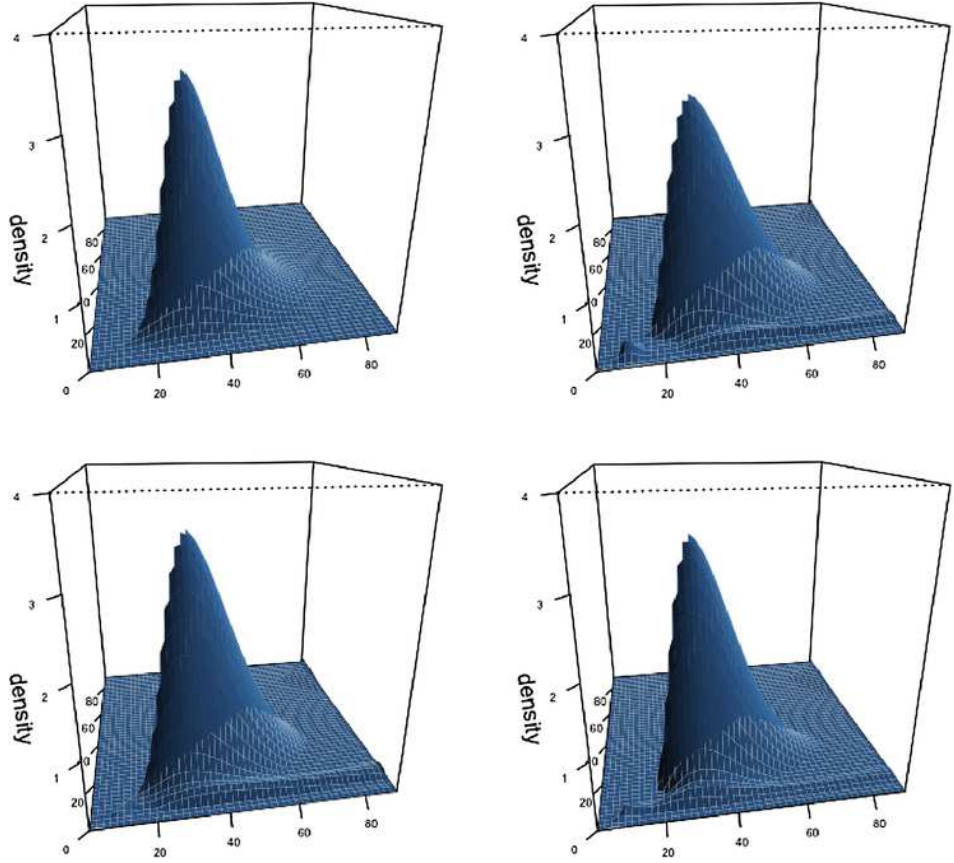}

\caption{Unimodal case: upper left displays the true density
of $W_{15}(2 I_2)$, upper right shows an estimate with $n = 500$, lower
left with $n = 1\mbox{,}000$ and lower right with $n = 2\mbox{,}000$.}\label{fig1}
\end{figure}

\begin{figure}

\includegraphics{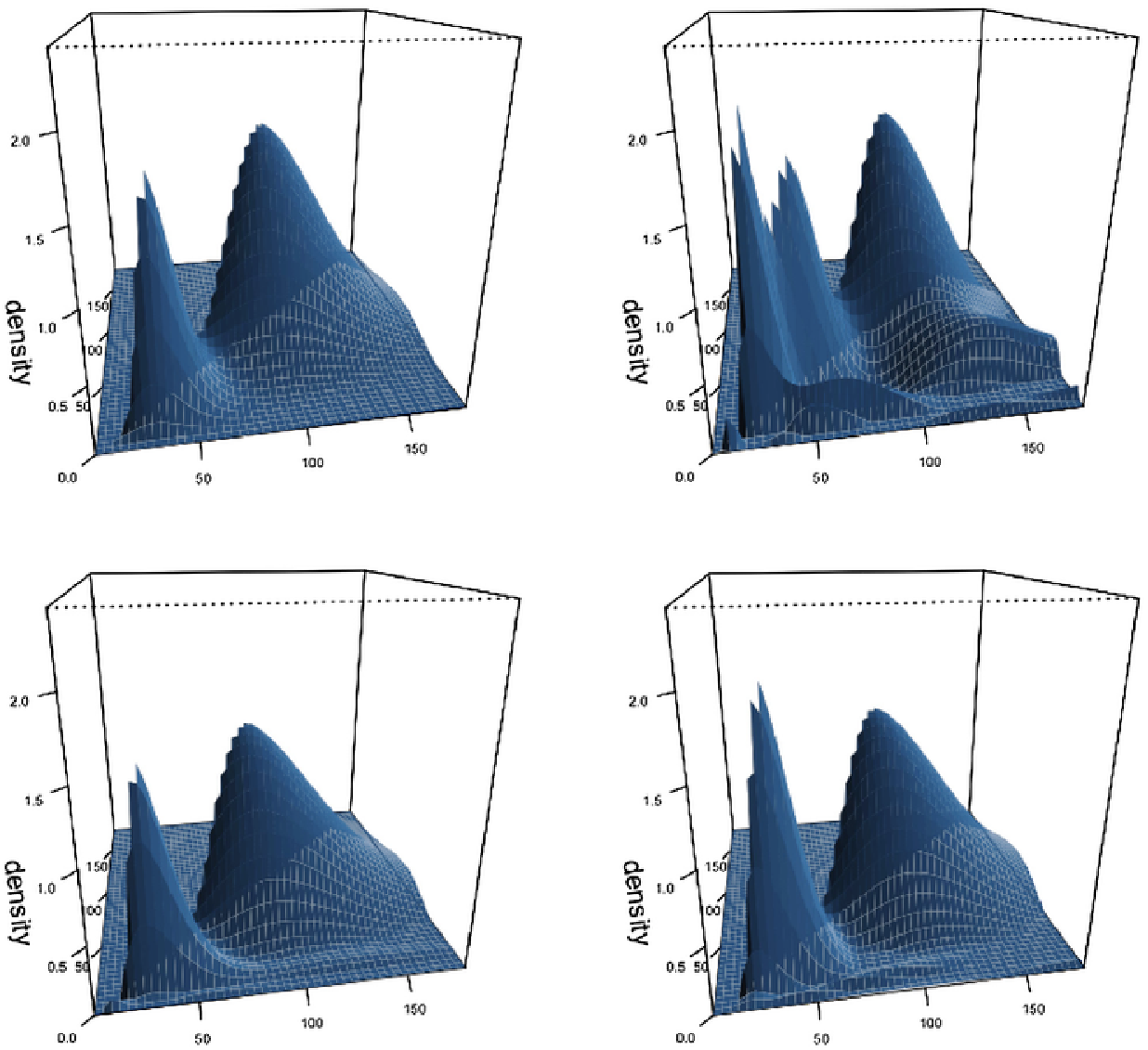}

\caption{Bimodal case: upper left displays the true density
$0.5 W_{15}(2 I_2) + 0.5 W_{15}(6 I_2)$, upper right shows an estimate
with $n = 500$, lower left with $n = 1\mbox{,}000$ and lower right with $n =
2\mbox{,}000$.}\label{fig2}
\end{figure}

As examples, we consider a unimodal mixing density $W_{15}(2 I_2)$
and a~bimodal density $ 0.5 W_{15}(2 I_2) + 0.5 W_{15}(6 I_2)$.
Figure \ref{fig1} show the results for the unimodal case
whereas Figure \ref{fig2} show the results for the bimodal case. In
each of these plots the
%
%
domain consists of the two eigenvalues starting with the largest.
One can see that the general shapes of the estimators become closer to
that of the true density
as $n$ increases.

\subsection{Application to stochastic volatility}
\label{subsecstock}
Stochastic volatility using the Wishart distribution is of much interest
in finance; see, for example,~\cite{amu} and~\cite{gs}.
In particular, this entails a situation precisely of the form (\ref{marginal}).
Let us apply this to the situation where we are interested in
estimating the mixing density.

Although our methods can be applied to a portfolio of many assets, let
us restrict
ourselves to two assets since this would be the smallest multivariate
example. Indeed,
let $S^1_j$ and $S^2_j$ denote the daily\vspace*{1pt} closing stock prices of
Samsung Electronics (005930.KS) and LG-display (034220.KS),
respectively, traded on the Korea Stock Exchange (KSC) for 2010,
where the data can be easily accessed on public financial websites.
We will\vspace*{1pt} assume as usual that $Q^k_j = \log(S^k_{j+1}/S^k_j
)$
follows a bi-variate normal distribution for $k = 1, 2$.
We transform the daily data to weekly data and compute the weekly
$2\times2$
covariance matrix $Y_i$
for $i=1,\ldots,52$. In case a week has
a holiday, we repeat the last previous observation. Under the usual assumptions
this would constitute observations from a mixture model (\ref
{marginal}) with a standard Wishart distribution
with four degrees of freedom.

%
\begin{figure}

\includegraphics{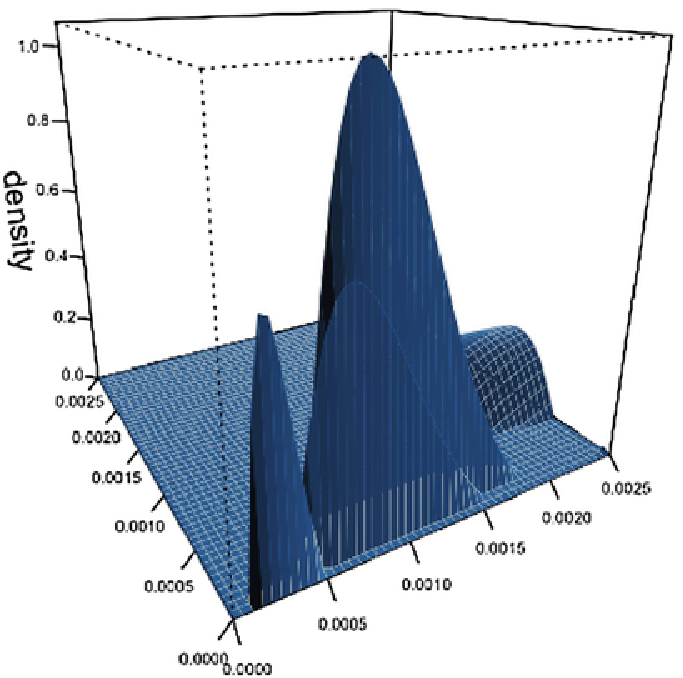}

\caption{A density estimator for the weekly covariance matrix of
stock prices.}\label{fig3}
\end{figure}

Figure \ref{fig3} plots the mixing density estimator corresponding to
the two eigenvalues. One can see that there are two peaks, suggesting a
possible bimodal stochastic volatility mixing density in the
eigenvalues.

\section{Proof of upper bound}
\label{secproof-upper}

The strategy here is, first, to decompose the integrated
mean-squared error into its variance and bias components,
%
%
\begin{eqnarray}
\label{ims}
\mE\|f_n - f\|^2 & = &
\mE\|(f_n - \mE f_n) + (\mE f_n - f)\|^2 \nonumber\\[-8pt]\\[-8pt]
& = & \mE\|f_n - \mE f_n\|^2 + \|\mE f_n - f\|^2,
\nonumber
\end{eqnarray}
and, last, to estimate each component separately using estimates
based on the Plancherel formula and the inversion formula for the
Helgason--Fourier transform.

\subsection{The integrated bias}
\label{ssecbias}

%
\begin{lemma}
\label{lemma-bias}
Suppose that $f \in\F(Q)$ and
$\varphi> \dim\Pm/2$. Then
\[
\|\mE f_n - f\|^2 \ll T^{-\varphi}.
\]
\end{lemma}
\begin{pf}
We have for $x \in\Pm$
%
%
\begin{eqnarray}
\label{bias-1}\quad
\mE f_n(x) - f(x)
&=& \int_{\C^m(\rho,T)} \int_{\bar{k} \in K/M}
\He f(s,\bar{k})
p_s(\bar{k}'x\bar{k}) \,{\mathrm{d}_*}\bar{k} \,{\mathrm{d}_*}s
\nonumber\\
&&{} - \int_{\C^m(\rho)} \int_{\bar{k} \in K/M}
\He f(s,\bar{k})
p_s(\bar{k}'x\bar{k}) \,{\mathrm{d}_*}\bar{k} \,{\mathrm{d}_*}s
\\
&=&-\int_{\lambda_s > T, \re(s)=-\rho} \int_{\bar{k} \in K/M}
\He f(s,\bar{k})
p_s(\bar{k}'x\bar{k}) \,{\mathrm{d}_*}\bar{k} \,{\mathrm{d}_*}s .\nonumber
\end{eqnarray}
Applying the Plancherel formula, we obtain
\[
\|\mE f_n - f\|^2 = \int_{\lambda_s > T,\re(s)=-\rho} \int_{K/M}
|\He f(s,\bar{k})|^2 \,{\mathrm{d}_*}\bar{k} \,{\mathrm{d}_*}s.
\]
Consequently,
%
%
\begin{eqnarray}
\label{i-bias}
\|\mE f_n - f\|^2 & = &
\int_{\lambda_s \geq T, \re(s)=-\rho} \int_{K/M}
|\He f(s,\bar{k})|^2 \,{\mathrm{d}_*}\bar{k} \,{\mathrm{d}_*}s \nonumber\\
&\le& T^{-\varphi}
\int_{\lambda_s \geq T, \re(s)=-\rho} \int_{K/M}
\lambda_s^{\varphi}
|\He f(s,\bar{k})|^2 \,{\mathrm{d}_*}\bar{k} \,{\mathrm{d}_*}s \\
&\le& T^{-\varphi}
\int_{\C^m(\rho)} \int_{K/M} \lambda_s^{\varphi}
|\He f(s,\bar{k})|^2\, {\mathrm{d}_*}\bar{k} \,{\mathrm{d}_*}s,
\nonumber
\end{eqnarray}
where we use the fact that
\[
\lambda_s^{\varphi} |\He f(s,k)|^2 \equiv
|\lambda_s^{\varphi/2} \He f(s,k)|^2 =
|\He(\Delta^{\varphi/2} f)(s,k)|^2.
\]
Therefore
\begin{eqnarray*}
\|\mE f_n - f\|^2 &\le& T^{-\varphi}
\int_{\C^m(\rho)} \int_{\bar{k}\in K/M}
|\He(\Delta^{\varphi/2} f)(s,k)|^2 \,{\mathrm{d}_*}\bar{k} \,{
\mathrm{d}_*}s \\
& = & T^{-\varphi} \int_{\Pm} |\Delta^{\varphi/2}f(w)|^2
\,{\mathrm{d}_*}w,
\end{eqnarray*}
where the equality follows from the Plancherel formula.
By assumption, $f \in\F(Q)$,
the latter integral is bounded above by $Q$, so we obtain
\[
\|\mE f_n - f\|^2 \leq Q T^{-\varphi},
\]
and the proof is complete.
\end{pf}

\subsection{The integrated variance}
\label{ssecvar}

To obtain bounds for the integrated variance, several preliminary
calculations are needed. In particular, we begin with the variance
calculation of the empirical Helgason--Fourier transform, which has
similarities to the usual empirical characteristic function.
%
%
\begin{lemma}
\label{lememp-helgason}
For $s \in\C^m(\rho)$ and $k \in K/M$,
\[
\mE|\He_n r(s,k) - \mE\He_n r(s,k)|^2 =
\frac1n \bigl(|\He r(-2\rho,k)|^2-|\He r(s,k)|^2\bigr).
\]
\end{lemma}
\begin{pf}
By (\ref{emp}),
%
%
\begin{eqnarray}
\label{var-1}
|\He_n r(s,k)|^2 &=& \overline{\He_n r(s,k)} \He_n r(s,k)
\nonumber\\
&=& \frac1{n^2} \sum_{j,\ell=1}^n \overline{p_s(k'Y_jk)}
p_s(kY_\ell
k) \\
&=& \frac1{n^2} \Biggl\{\sum_{j=1}^n |p_s(k'Y_jk)|^2
+ \sum_{j\neq\ell}\overline{p_s(k'Y_jk)}p_s(k'Y_\ell
k)\Biggr\}.\nonumber
\end{eqnarray}
Observe also that
\begin{eqnarray*}
p_s(w)\overline{p_s(w)} & = &
|w_1|^{s_1}\cdots|w_m|^{s_m}\overline{|w_1|^{s_1}\cdots|w_m|^{s_m}}
\\
& = & |w_1|^{2 \re(s_1)}\cdots|w_m|^{2 \re(s_m)} \\
& = & p_{-2\rho}(w)
\end{eqnarray*}
since $\re(s)=-\rho$. Applying this result to (\ref{var-1}) and
taking expectations, we obtain
\begin{eqnarray*}
\mE|\He_n r(s,k)|^2
&=& \frac1{n^2}\mE\Biggl\{\sum_{j=1}^n|p_s(k'Y_jk)|^2
+ \sum_{j\neq\ell}\overline{p_s(k'Y_jk)}p_s(k'Y_\ell k)\Biggr\} \\
&=& \frac1{n^2}\Biggl\{ \sum_{j=1}^n \mE p_{-2\rho}(k'Y_jk) +
\sum_{j\neq\ell} \overline{\mE p_s(k'Y_jk)}
\mE p_s(k'Y_\ell k)\Biggr\} \\
&=& \frac1{n}\He r(-2\rho,k) +
\frac{n-1}n |\He r(s,k)|^2,
\end{eqnarray*}
where the last equality follows from the fact that $Y_1,\ldots,Y_n$
are independent and identically distributed as $Y$, and because
$\mE p_s(k'Yk) = \He r(s,k)$.
\end{pf}

Following Terras \cite{Terras-II}, pages 34 and 35, let
\[
A^+ = \{a = \diag(a_1,\ldots,a_m) \in A\dvtx a_1 > \cdots> a_m \}
\]
denote the positive Weyl chamber in $A$.
For $a = \diag(a_1, \ldots, a_m) \in A$,
let $\mathrm{d} a = \prod_{j=1}^m a_j^{-1} d a_j$ and set
\[
\gamma(a) = \prod_{j=1}^m a_j^{-(m-1)/2} \prod_{1 \le i < j \le m} |a_i
- a_j|,
\]
and define the normalizing constant $b_m$ by
$b_m^{-1} = \pi^{-(m^2+m)/4} \prod_{j=1}^m j \Gamma(j/2)$.
Denote ${\mathrm{d}_*}a = b_m \gamma(a) \,\mathrm{d} a$.
%
%
\begin{lemma}
\label{lemma-var}
As $T \to\infty$,
\[
\mE\|f_n - \mE f_n\|^2 \ll{\sup_{s \in\C^m(\rho)}}
|
\what_{Z}(s)|^{-2} \frac{T^{\dim\Pm/2}}{n}.
\]
\end{lemma}
\begin{pf} By the Plancherel formula,
\begin{eqnarray*}
&&\mE\|f_n - \mE f_n\|^2 \\
&&\qquad= \int_{\C^m(\rho,T)} \int_{K/M}
\mE|\He_n f(s,k) - \mE\He_n f(s,k)|^2
\,\mathrm{d}\bar{k} \,{\mathrm{d}_*}s \\
&&\qquad= \int_{\C^m(\rho,T)} \int_{K/M}
\mE|\He_n r(s,k) - \mE\He_n r(s,k)|^2
\,\mathrm{d}\bar{k} |\what_{Z}(s)|^{-2} \,{\mathrm{d}_*}s \\
&&\qquad\leq
\frac1n \sup_{\lambda_s < T,
\re(s)=-\rho}|\what_{Z}(s)|^{-2} \int_{K/M}
|\He r(-2\rho,\bar{k})|^2\,\mathrm{d}\bar{k}\int_{\C^m(\rho,T)}
{\mathrm{d}_*}s \\
&&\qquad\ll{\sup_{\lambda_s < T,\re(s)=-\rho}}
|\what_{Z}(s)|^{-2} \frac{T^{\dim\Pm/2}} n
\end{eqnarray*}
as $T \to\infty$.

Choose $a \in A^+$.
Observe that
%
%
\begin{equation} \label{lemhbound}
p_{-2\rho}(a)
= a_1^{-(m-1)+(m-1)/2} \cdots a_{m-1}^{-1+(m-1)/2} a_m^{(m-1)/2}
\le1.
\end{equation}
Since $p_{-2\rho}(k'ak)$ is a continuous function of $k$ on a compact
set $K$,
$p_{-2\rho}(k'ak)$ is uniformly bounded on $K$ such\vadjust{\goodbreak}
$|p_{-2\rho} (k'ak)| \le C$ on $K$.
Since $f$ is a density so that $r$ is also a density, we have
\begin{eqnarray*}
|\He r(-2\rho, \I_m)|
&=& \biggl|\int_\Pm r(y) \overline{p_{-2\rho} (y)}\, {\mathrm{d}_*}y
\biggr| \\
&=& \int_{A^+} \int_K r(k' a k) \overline{p_{-2\rho} (k' a k)}
\,{\mathrm{d}_*}a\,
{\mathrm{d}_*}k \\
&\le& \int_{A^+} \int_K r(k' a k) |\overline{p_{-2\rho} (k' a
k)}| \,{\mathrm{d}_*}a\, {\mathrm{d}_*}k \\
&\le& C \int_{A^+} \int_K r(k' a k) \,{\mathrm{d}_*}a\, {\mathrm{d}_*}k \\
&=& C \int_\Pm r(y)\, {\mathrm{d}_*}y \\
&=& C.
\end{eqnarray*}
Hence, it follows from continuity and the compactness of $K/M$
\[
\int_{K/M} |\He r(-2\rho,\bar{k})|^2\,\mathrm{d}\bar{k} < \infty,
\]
which has been used in the above calculation.

In addition, we use the fact that, as $T \to\infty$,
\[
{\sup_{\C^m(\rho,T)}} |c_m(s)|^{-2} \ll T^{m(m-1)/4},
\]
a result which follows from Proposition 7.2 of Helgason \cite{he2},
page 450.
\end{pf}

The proof of the upper bound can now be obtained by applying Lem\-mas~\ref{lemma-bias} and \ref{lemma-var} to (\ref{ims}) and setting
$T \asymp(\log n)^2$.

\section{Proof of lower bound}
\label{secproofs-lower}

We need to provide some detailed calculations, and the
essence of the proof is contained for the case $m=2$;
hence we will keep this assumption for the remainder of this paper.
The generalization to $m > 2$ may
be obtained by using higher order hyperbolic spherical coordinates.
In this section, we assume that $\psi$ is a $K$-invariant function
defined on $\P$.

\subsection{Convolution and Helagson--Fourier transform in polar coordinate}

For $y \in\P$, let $y = k' a k $ with $a = \diag(a_1, a_2) \in A^+$,
$k \in K$
so that $\psi(y) = \psi(a)$.
Let
\[
a = D_{u_1} e^{u_2}
\]
with $D_z = \diag(e^{z}, e^{-z} )$ for $z \in\R$,
and write
\[
\psi(u) = \psi(D_{u_1} e^{u_2}).
\]
By a change of variables,
%
%
\begin{equation}
\label{eqy-to-polar}
\int_{\Pos_2} \psi(y) \,{\mathrm{d}_*}y
= \int_{\D} \psi(u)\, {\mathrm{d}_*}u,
\end{equation}
where
\[
\D= \{ u\dvtx u_1 \in\R^+ \mbox{ and } u_2 \in\R\}
\]
and
\[
{\mathrm{d}_*}u = 4 \pi\sinh u_1 \,\mathrm{d}u_1 \,\mathrm{d}u_2.
\]

Denote
$
k_\theta=
\bigl[{\cos\theta\atop\sin\theta} \enskip{{-}\sin\theta
\atop\cos\theta}\bigr]
$.
The next lemma is straightforward so we shall omit the proof.
%
%
\begin{lemma} \label{lemMatrixEquation}
For $u_1 \in\R^+$, $v_1 \in\R^+$ and $\theta\in[0, 2\pi]$,
the matrix equation
\[
k_\xi D_R k_\xi' = D_{-u_1/2} k_\theta D_{v_1} k_\theta' D_{-u_1/2}
\]
has a solution $R^* = R^*(u_1, v_1, \theta)$ and
$\xi^* = \xi^*(u_1, v_1, \theta)$. Further, $\cosh R^*$ has
minimum and maximum values $\cosh(u_1-v_1)$ and
$\cosh(u_1+v_1)$, respectively, and $R^*$ and $\xi^*$ can be
defined uniquely.
\end{lemma}

In general, if both $f$ and $g$ are $K$-invariant functions on $\Pm$,
then $f * g$ is also $K$-invariant and $f * g = g * f$.
Hence, $\psi* w$ is $K$-invariant and $\psi*w = w*\psi$
due to $K$-invariance of $\psi$ and $w$.
From this and Lemma \ref{lemMatrixEquation},
we can define
\[
\Psi_{u_1 v_1}(z) =
\frac1{2\pi} \int_0^{2\pi}
\psi(R^*(\theta, u_1, v_1), z ) \,\mathrm{d}\theta
\]
for $z \in\R$, $u_1 \in\R^+$, $v_1 \in\R^+$.
Denoting by $W$ the distribution function corresponding to
the standard Wishart density $w(u)$ with respect to the measure~${\mathrm{d}_*}u$,
we have that for $v \in D$,
%
%
\begin{equation}
\label{eqconv-polar}
(\psi* w)(v) = \int_{\D} \Psi_{u_1 v_1}(v_2 - u_2) \,\mathrm{d}W(u).
\end{equation}

The Laplacian for $K$-invariant functions in polar coordinate is given by
$\Delta= \Delta_{u_1} + \Delta_{u_2}$,
where
\[
\Delta_{u_1} = -\coth(u_1) \frac\partial{\partial u_1} - \frac
{\partial^2} {\partial u_1^2},\qquad
\Delta_{u_2} = -\frac{\partial^2} {\partial u_2^2},
\]
and the spherical function is given by
%
%
\begin{equation}
\label{eqpolar-hs}
h_s(u) = P_{s_1} (\cosh u_1) e^{(s_1+2s_2) u_2}\qquad
\mbox{for } u \in\D
\end{equation}
with $P_s$, the Legendre function; see Terras \cite{Terras-II}.
It can be seen that
\[
\Delta_{u_1} h_s = - \tfrac{1}{2} s(s+1) h_s,\qquad
\Delta_{u_2} h_s = - \tfrac{1}{2}(s_1+2s_2)^2 h_s,
\]
so that $\Delta h_s = \lambda_s h_s$
with $\lambda_s = - \half\{ (s_1+2s_2)^2 + s(s+1) \}$.\vadjust{\goodbreak}

The Helgason--Fourier transform of $\psi$ is given by
%
%
\begin{equation}
\label{eqpolar-HF}
\hat\psi(s) = \int_\D\psi(u) P_{\bar s_1}( \cosh u_1)
e^{(\overline{s_1+2s_2})u_2} \,{\mathrm{d}_*}u
\end{equation}
from (\ref{eqy-to-polar}) and (\ref{eqpolar-hs}).
Suppose $\psi$ is separable so that $\psi(u) = \psi_1(u_1) \psi_2(u_2)$.
Then, (\ref{eqpolar-HF}) implies
%
%
\begin{equation}
\label{eqHF-sep}
\hat\psi(s) = {\mathcal M} \psi_1(\bar s_1) \L\psi
_2(s_1+2s_2) ,
\end{equation}
where $\L$ and $\M$ denote, respectively, the Laplace transform and the
Mehler--Fock transform;
see Terras \cite{Terras-II}.

\subsection{\texorpdfstring{$\chi^2$-divergence}{chi^2-divergence}}

Choose $\psi_1$ as the perturbation function as in Fan~\cite{fan}.
Then, one can construct $\psi(u)$ satisfying the following conditions:
\begin{longlist}[(P5)]
\item[(P1)]
$\psi$ is $K$-invariant and separable with $\psi(u) = \psi_1(u_1)
\psi_2(u_2)$
for $u \in D$.
\item[(P2)]
$\L\psi_2(it) = 0$ for $t \notin[1,2]$.
\item[(P3)]
$\psi\in\F(Q)$.
\item[(P4)]
$\psi_1(u_1) = O(\cosh^{-m_0} u_1)$ and $\psi_2(u_2) = O
(e^{-m_0 |u_2|})$,
where $0 < \xi< 1$ and $m_0 \xi> (N-1)(1-\xi)/2$.
\item[(P5)]
$\int_\D\psi(u) \,{\mathrm{d}_*}u = 0$.
\end{longlist}

Let $p_b$ be a density on $\R$ such that
$p_b$ is sufficiently smooth and satisfies $p_b(u_2) = c_b \exp(- b |u_2|)$,
$|u_2| \ge c_0$, where $c_b$ is a normalizing constant.
Define the function
\[
f^0(u) = C_b (\cosh u_1)^{-b} p_b(u_2) ,
\]
where $C_b = c_b (b-1)/(2\pi)$ for $b > 1$.

For a function $g\dvtx \P\rightarrow\R$ and $\delta> 0$, define
\[
g^\delta(y) = g\bigl( |y|^{(\delta-1)/2} y \bigr) \qquad\mbox{for } y \in\P
\]
so that
\[
\psi(u) = \psi(u_1, \delta u_2) \qquad\mbox{for } u \in\D.
\]
%
%
\begin{proposition}
\label{propchisq-1}
Suppose \textup{(P1)--(P5)} hold.
For a pair of densities
\[
f^0 \quad\mbox{and}\quad f^n = f^0 + C_\psi\delta^{-\varphi+ 1/2} \psi^\delta,
\]
the $\chi^2$-divergence between $g^0 = f^0 * w$ and $g^n = f^n * w$
satisfies
\[
\chi^2(g^0,g^n) \le C/n
\]
provided that $b < \frac1 2 \min( 3\pi, (N-1)(1-\xi) - 1
)$.
\end{proposition}

The fact that $(\log n)^{-2\varphi}$ is a lower bound follows from
Proposition \ref{propchisq-1} whose proof follows from a sequence of
lemmas below.

\subsection{Perturbing function}

Denote $(q_1, q_2) = (s_1, s_1 + 2s_2)$
and $\beta_j = \im(q_j)$ for $j = 1, 2$.
Note that $q_2 = i \beta_2$ for $s \in C^2(\rho)$.
%
%
\begin{lemma} \label{lemPsi}
Suppose \textup{(P1)} holds.
Then
$\| \psi^\delta\| = \delta^{-1} \| \psi\|$.
\end{lemma}
\begin{pf}
By (\ref{eqHF-sep}) and the change of variable $u_2 \mapsto\delta
u_2$, we obtain
\[
\hat{\psi}^\delta(s)
= \M\psi_1(\bar q_1) \{\delta^{-1} \L\psi_2(\bar
{q}_2/\delta)\}.
\]
The desired result follows from
the Plancherel formula (\ref{hf-plancherel}) and change of variable $s
\mapsto q$.
\end{pf}
%
%
\begin{lemma} \label{lemPsiF}
Suppose \textup{(P1), (P2)} and \textup{(P3)} hold.
Then there exists a~posi\-tive constant $C_\psi$ such that
$C_\psi\delta^{-\varphi+ 1/2} \psi^\delta\in\F(Q)$.
\end{lemma}
\begin{pf}
Suppose $s \in C^2(\rho) \cap{\mathcal S}$.
Observe that
\[
\lambda_s = -\tfrac12\{ s_1(s_1+1) + (s_1+2s_2)^2 \}
= \tfrac12 \bigl( \beta_1^2 + \beta_2^2 + \tfrac14 \bigr)
\]
and that for $\delta\ge1$,
\[
\tfrac12 \bigl( \beta_1^2 + (\delta\beta_2)^2 + \tfrac14 \bigr)
\le\delta^2 \lambda_s.
\]
Now, the Plancherel formula (\ref{hf-plancherel}) and change of
variable $s \mapsto q$ gives
\[
\| \Delta^{\varphi/2} \psi^\delta\|^2
\le C \delta^{2\varphi- 1} \| \Delta^{\varphi/2} \psi\|^2.
\]
A suitable choice of $C_\psi$ gives the desired result.
\end{pf}
%
%
\begin{lemma} \label{lemPsiC}
Under \textup{(P1)} and \textup{(P2)},
$\| \psi^\delta* w \|^2 \le C \delta e^{-3 \pi\delta} \| \psi\|^2$.
\end{lemma}
\begin{pf}
For $s \in C^2(\rho) \cap{\mathcal S}$,
\[
|\hat w(s)|^2 \le C e^{- \pi(\alpha_1 + 2 \alpha_2)}
= C e^{-2\pi\beta_2}
\le C e^{-3\pi}.
\]
The desired result follows from the inequality $e^{- \pi(s_1 + 2 s_2)
\delta} \le C e^{-3\pi\delta}$, the Plancherel formula
(\ref{hf-plancherel}) and change of variable $s \mapsto q$.
\end{pf}

\subsection{Tail behavior}

%
\begin{lemma}
\label{lemw-tail}
For $c_1 \ge0$ and $c_2, c_3 \in\R$,
\[
\int_{\{ u_1 > {c_1}, c_2 < u_2 < c_3 \}} \mathrm{d} W(u)
= 4\pi\int_{\{c_2 < u_2 < c_3\}} \exp\{ (N-1)u_2 - e^{u_2}
\cosh
{c_1} \} \,\mathrm{d} u_2.
\]
\end{lemma}
\begin{pf}
Change of variable gives the desired result.\vadjust{\goodbreak}
\end{pf}

If $U$ has the distribution function $W$,
then $(U_1, \delta U_2)$ has the distribution function $W_\delta$
with density $\delta^{-1} w(u_1, u_2/\delta)$.
By (\ref{eqconv-polar}) and a change of variables,
%
%
\begin{equation}
\label{eqconv-polar-delta}
(\psi^\delta* w)(v) = \int_{\D}\Psi_{u_1v_1}(\delta v_2 - u_2)
\,\mathrm{d}W_\delta(u).
\end{equation}
%
%
\begin{lemma}
\label{lemg0}
We have $g^0 = f^0 * w$ is $K$-invariant, and
as $v_1 \to\infty$ and $|v_2| \to\infty$
\[
g^0(v) \ge C e^{-b(v_1+|v_2|)}.
\]
\end{lemma}
\begin{pf}
Choose a constant $C_{1/2}$
such that
\[
\int_{\{0 \le u_1 < C_{1/2}, |u_2| < C_{1/2} \}} \mathrm{d} W(u) \ge1/2
\]
and $p_b(u_2)\,{=}\,c_b \exp({-}b |u_2|)$ for $|u_2|\,{\ge}\,C_{1/2}$.
The desired result follows from~(\ref{eqconv-polar}) and Lemma \ref
{lemMatrixEquation}.
\end{pf}
%
%
\begin{lemma}
\label{lembound-munu}
Define
\[
\eta_\delta(v) = \int_\D\cosh^{-m_0} (u_1-v_1) e^{-m_0 |v_2 -
u_2|}
\,\mathrm{d} W_\delta(u).
\]
Then, there exist $M$ and a constant $C$ such that
when $v_1 \ge M$ and \mbox{$|v_2|/\delta\ge M$},
\[
\eta_\delta(v)
\le C \exp\bigl\{-\tfrac12 (N-1)(1-\xi) (v_1 + |v_2|/\delta)
\bigr\}
\]
for all $\delta\ge1$ provided that $0 < \xi< 1$ and $m_0 \xi>
(N-1)(1-\xi)/2$.
\end{lemma}
\begin{pf}
Denote
\begin{eqnarray*}
J_{11} &=& \{u \dvtx |u_1 - v_1| \le\xi v_1, |v_2 - \delta u_2| \le\xi
|v_2| \},\\
J_{12} &=& \{u \dvtx |u_1 - v_1| \le\xi v_1, |v_2 - \delta u_2| > \xi
|v_2| \},\\
J_{21} &=& \{u \dvtx |u_1 - v_1| > \xi v_1, |v_2 - \delta u_2| \le\xi
|v_2| \},\\
J_{22} &=& \{u \dvtx |u_1 - v_1| > \xi v_1, |v_2 - \delta u_2| > \xi|v_2|
\}
\end{eqnarray*}
and for $i,j = 1,2$,
\[
I_{ij} = \int_{J_{ij}} \cosh^{-m_0} (u_1-v_1) e^{-m_0 |v_2 - u_2|}\, \mathrm{d}
W_\delta(u).
\]

Consider $I_{11}$. Suppose $v_2 < 0$.
If $|u_1 - v_1| \le\xi v_1$ and $|v_2 - \delta u_2| \le\xi|v_2|$,
then
\[
(1-\xi) (v_1+|v_2|/\delta) \le u_1 - u_2 \le(1+\xi) (v_1 +
|v_2|/\delta).
\]
Let $\zeta= (1-\xi) (v_1+|v_2|/\delta)$ and $\lambda= e^\zeta$.
Note that
\[
I_{11}
\le
C \biggl(\int_{\{u_1 \ge0, u_2 < -\zeta\}} \mathrm{d} W(u)
+
\int_{\{ u_1 \ge u_2 + \zeta, -\zeta\le u_2 <0 \}} \mathrm{d}
W(u)\biggr)
:= I_1^- + I_2^-.\vadjust{\goodbreak}
\]
Since
$
e^{u_2} \cosh(u_2 + \zeta)
= \frac12 ( \lambda e^{2u_2} + \lambda^{-1}),
$
then we have
\[
I_1^- = 4\pi\int_{\{ 0 \le t < 1/\lambda\}} t^{N-2} e^{-t}\,\mathrm{d}t
= O \bigl( \lambda^{-(N-1)} \bigr)
\]
and by Lemma \ref{lemw-tail} with $c_1 = u_2 + \zeta$, $c_2 = -\zeta$
and $c_3 = 0$
\[
I_2^-
= C \int_{\{ -\zeta\le u_2 <0 \}}
e^{(N-1)u_2} \exp\bigl( -e^{u_2} \cosh(u_2 + \zeta) \bigr) \,\mathrm{d} u_2
= O \bigl( \lambda^{-(N-1)/2} \bigr).
\]
Hence, if $v_2\,{<}\,0$, then
$
I_{11}\,{=}\,O ( \lambda^{-(N-1)/2} )
$
as $v_1\,{+}\,|v_2|/\delta\,{\to}\,\infty$.
Suppose \mbox{$v_2\,{\ge}\,0$}.
If $|v_2 - \delta u_2| \le\xi|v_2|$ and $|u_1 - v_1| \le\xi v_1$,
then
\[
(1-\xi) (v_1 + |v_2|/\delta) \le u_1 + u_2 \le(1+\xi) (v_1 +
|v_2|/\delta).
\]
Observe that
\[
I_{11}
\le C \biggl(\int_{\{u_1 \ge0, u_2 \ge\zeta\}} \mathrm{d}W(u)
+
\int_{\{u_1 \ge\zeta- u_2, 0 \le u_2 < \zeta\}}\mathrm{d}W(u)\biggr)\\
:= I_1^+ + I_2^+.
\]
Since $e^{u_2/2} \cosh(\zeta- u_2) = \frac12 ( \lambda+
\lambda
^{-1} e^{2u_2} )$,
then we have
\[
I_1^+ = 4\pi\int_{\{ t \ge\lambda\}} t^{N-2} e^{-t}\,\mathrm{d}t = O(
e^{-\lambda})
\]
and by Lemma \ref{lemw-tail}, with $c_1 = u_2 - \zeta$, $c_2 = 0$ and
$c_3 = \zeta$,
\[
I_2^+
= C \int_{\{ 0 \le u_2 < \zeta\}} e^{(N-1)u_2} \exp\bigl( -e^{u_2}
\cosh(\zeta- u_2) \bigr) \,\mathrm{d} u_2
= O(e^{-\lambda}).
\]
Hence, for $v_2 \in\R$,
\[
I_{11} = O \bigl( \lambda^{-(N-1)/2} \bigr)
= O \bigl( \exp\bigl\{ - \tfrac12 (N-1)(1-\xi) (v_1+|v_2|/\delta)
\bigr\} \bigr)
\]
as $v_1 + |v_2|/\delta\to\infty$.

Consider $I_{21}$.
Let $\zeta= (1-\xi) |v_2|/\delta$ and $\lambda= e^\zeta$.
Suppose $v_2 < 0$.
On $J_{21}$, $u_2 \le-(1-\xi)|v_2|/\delta= -\zeta$ on $J_{21}$.
Note that
$
\cosh(\xi v_1)^{-m_0} \le(\frac12 e^{\xi v_1})^{-m_0}
\le2^{-m_0} e^{-(1/2)(N-1)(1-\xi)v_1}.
$
If follows from these and Lemma \ref{lemw-tail} that
\[
I_{21}
= O \bigl( \exp\bigl\{ - \tfrac12 (N-1)(1-\xi) (v_1 + |v_2|/\delta)
\bigr\} \bigr),
\]
and if $v_2 \ge0$,
\[
I_{21}
= O \bigl( \exp\bigl\{ - \tfrac12 (N-1)(1-\xi) (v_1 + |v_2|/\delta)
\bigr\} \bigr).
\]

Observe that
$m_0 \xi|v_2| \ge m_0 \xi|v_2|/\delta\ge\frac1 2 (N-1)(1-\xi)
|v_2|/\delta$.
It follows from this and Lemma \ref{lemw-tail} that
\begin{eqnarray*}
I_{12}
&\le& C e^{-m_0 \xi|v_2|}
\int_{u_2 \in\R} \int_{u_1 > (1 - \xi) v_1} \mathrm{d}W(u)
\\
&=& O \biggl( \exp\biggl\{ - \frac12 (N-1)(1-\xi) (v_1 +
|v_2|/\delta)
\biggr\} \biggr)
\end{eqnarray*}
and
\begin{eqnarray*}
I_{22}
&\le& C e^{-m_0 \xi|v_2|} (\cosh(\xi v_1) )^{-m_0}
\int_{\D} \mathrm{d}W(u)
\\
&=& O \biggl( \exp\biggl\{-\frac12 (N-1)(1-\xi) (v_1 + |v_2|/\delta
)\biggr\} \biggr).
\end{eqnarray*}
This completes the proof of Lemma \ref{lembound-munu}.
\end{pf}

\subsection{\texorpdfstring{Proof of Proposition \protect\ref{propchisq-1}}{Proof of Proposition 7.2}}

From (P5), one can choose $C_\psi$ sufficiently small such that $f^n$
is a density.
By (\ref{eqy-to-polar}) and (\ref{eqconv-polar-delta}),
\begin{eqnarray*}
\chi^2(g^0,g^n)
&=& C_\psi^2 \delta^{-2\varphi+ 1}
\int_\D\frac{\{ (\psi^\delta* w)(v) \}^2}{g^0(v)}
\,{\mathrm{d}_*}v
\\
&=& C_\psi^2 \delta^{-2\varphi}
\int_\D\frac{\{ \int_{\D} \Psi_{u_1 v_1}(v_2 - u_2)
\,\mathrm{d}W_\delta(u)
\}^2}
{g^0(v_1, v_2/\delta)} \,{\mathrm{d}_*}v.
\end{eqnarray*}
Let
\begin{eqnarray*}
J_{11} &=& \{v\dvtx 0 \le v_1 \le\delta, |v_2|/\delta\le\delta\},\qquad
J_{12} = \{v\dvtx 0 \le v_1 \le\delta, |v_2|/\delta> \delta\},\\
J_{21} &=& \{v\dvtx v_1 > \delta,|v_2|/\delta\le\delta\},\qquad
J_{22} = \{v\dvtx v_1 > \delta, |v_2|/\delta> \delta\},
\end{eqnarray*}
and for $i,j = 1,2$,
\begin{eqnarray*}
I_{ij} &=& \int_{J_{ij}} \frac{\{ \int_{\D} \Psi_{u_1
v_1}(v_2 -
u_2) \,\mathrm{d}W_\delta(u) \}^2}
{g^0(v_1, v_2/\delta)} \,{\mathrm{d}_*}v.
\end{eqnarray*}
By (\ref{eqconv-polar-delta}),
Lemmas \ref{lemPsiC}, \ref{lemg0} and the Plancherel formula
(\ref
{hf-plancherel}),
we obtain
\begin{eqnarray*}
I_{11}
&\le& C e^{2b\delta}
\int_\D\biggl\{\int_{\D} \Psi_{u_1 v_1}(v_2 - u_2) \,\mathrm{d}W_\delta(u)
\biggr\}
^2 \,{\mathrm{d}_*}v
\\
&=& C e^{2b\delta} \delta\| \psi^\delta* w \|^2
\\
&\le& C \delta^2 e^{(2b-3\pi)\delta}.
\end{eqnarray*}
Lemma \ref{lemMatrixEquation} and (P4) imply
%
%
\begin{equation}
\label{eqPsieta}
\int_{\D} \Psi_{u_1 v_1}(v_2 - u_2) \,\mathrm{d}W_\delta(u) \le C \eta_\delta(v).
\end{equation}
Let $c_1 = (N-1)(1-\xi) - b$.
It follows from (\ref{eqPsieta}), Lemmas \ref{lemg0} and \ref
{lembound-munu} that
\[
I_{12} = O \bigl( \delta e^{-(c_1-b) \delta} \bigr),\qquad
I_{21} = O \bigl(e^{-(c_1-b-1) \delta} \bigr),\qquad
I_{22} = O \bigl(\delta e^{-(2c_1-1) \delta} \bigr).
\]

Now letting $\varepsilon= \min(3\pi- 2b, (N-1)(1-\xi) - 2b -
1
) > 0$
and combining the above bounds, we obtain
\[
\chi^2(g^0, g^n)
\le C_\psi^2 \delta^{-2\varphi} ( I_{11} + I_{12} + I_{21} +
I_{22} )
\le C \delta^{-2\varphi+ 2} e^{-\varepsilon\delta}.
\]
Choosing $\delta= \varepsilon/\log n$, we have the desired result.


%
%

\printaddresses


\begin{thebibliography}{27}

\bibitem{amu}
%
\begin{barticle}[mr]
\bauthor{\bsnm{Asai},~\bfnm{Manabu}\binits{M.}},
\bauthor{\bsnm{McAleer},~\bfnm{Michael}\binits{M.}} \AND
\bauthor{\bsnm{Yu},~\bfnm{Jun}\binits{J.}}
(\byear{2006}).
\btitle{Multivariate stochastic volatility: A~review}.
\bjournal{Econometric Rev.}
\bvolume{25}
\bpages{145--175}.
\bid{doi={10.1080/07474930600713564}, issn={0747-4938}, mr={2256285}}
\bptok{imsref}%
\end{barticle}
%
\endbibitem

\bibitem{butsy}
%
\begin{barticle}[mr]
\bauthor{\bsnm{Butucea},~\bfnm{C.}\binits{C.}} \AND
\bauthor{\bsnm{Tsybakov},~\bfnm{A.~B.}\binits{A.~B.}}
(\byear{2008}).
\btitle{Sharp optimality in density deconvolution with dominating bias,
I, II}.
\bjournal{Theory Probab. Appl.}
\bvolume{52}
\bpages{24--39, 237--249}.
\bptok{imsref}%
\end{barticle}
%
\endbibitem

\bibitem{chentan}
%
\begin{barticle}[mr]
\bauthor{\bsnm{Chen},~\bfnm{Jiahua}\binits{J.}} \AND
\bauthor{\bsnm{Tan},~\bfnm{Xianming}\binits{X.}}
(\byear{2009}).
\btitle{Inference for multivariate normal mixtures}.
\bjournal{J. Multivariate Anal.}
\bvolume{100}
\bpages{1367--1383}.
\bid{doi={10.1016/j.jmva.2008.12.005}, issn={0047-259X}, mr={2514135}}
\bptok{imsref}%
\end{barticle}
%
\endbibitem

\bibitem{diggle}
%
\begin{barticle}[mr]
\bauthor{\bsnm{Diggle},~\bfnm{Peter~J.}\binits{P.~J.}} \AND
\bauthor{\bsnm{Hall},~\bfnm{Peter}\binits{P.}}
(\byear{1993}).
\btitle{A {F}ourier approach to nonparametric deconvolution of a density
estimate}.
\bjournal{J. Roy. Statist. Soc. Ser. B}
\bvolume{55}
\bpages{523--531}.
\bid{issn={0035-9246}, mr={1224414}}
\bptok{imsref}%
\end{barticle}
%
\endbibitem

\bibitem{fan}
%
\begin{barticle}[mr]
\bauthor{\bsnm{Fan},~\bfnm{Jianqing}\binits{J.}}
(\byear{1991}).
\btitle{On the optimal rates of convergence for nonparametric deconvolution
problems}.
\bjournal{Ann. Statist.}
\bvolume{19}
\bpages{1257--1272}.
\bid{doi={10.1214/aos/1176348248}, issn={0090-5364}, mr={1126324}}
\bptok{imsref}%
\end{barticle}
%
\endbibitem

\bibitem{gs}
%
\begin{barticle}[mr]
\bauthor{\bsnm{Gourieroux},~\bfnm{Christian}\binits{C.}} \AND
\bauthor{\bsnm{Sufana},~\bfnm{Razvan}\binits{R.}}
(\byear{2010}).
\btitle{Derivative pricing with {W}ishart multivariate stochastic volatility}.
\bjournal{J. Bus. Econom. Statist.}
\bvolume{28}
\bpages{438--451}.
\bid{doi={10.1198/jbes.2009.08105}, issn={0735-0015}, mr={2723611}}
\bptok{imsref}%
\end{barticle}
%
\endbibitem

\bibitem{haff1}
%
\begin{barticle}[mr]
\bauthor{\bsnm{Haff},~\bfnm{L.~R.}\binits{L.~R.}}
(\byear{1979}).
\btitle{Estimation of the inverse covariance matrix: Random mixtures of the
inverse {W}ishart matrix and the identity}.
\bjournal{Ann. Statist.}
\bvolume{7}
\bpages{1264--1276}.
\bid{issn={0090-5364}, mr={0550149}}
\bptok{imsref}%
\end{barticle}
%
\endbibitem

\bibitem{haff2}
%
\begin{barticle}[mr]
\bauthor{\bsnm{Haff},~\bfnm{L.~R.}\binits{L.~R.}}
(\byear{1980}).
\btitle{Empirical {B}ayes estimation of the multivariate normal covariance
matrix}.
\bjournal{Ann. Statist.}
\bvolume{8}
\bpages{586--597}.
\bid{issn={0090-5364}, mr={0568722}}
\bptok{imsref}%
\end{barticle}
%
\endbibitem

\bibitem{haff3}
%
\begin{barticle}[mr]
\bauthor{\bsnm{Haff},~\bfnm{L.~R.}\binits{L.~R.}}
(\byear{1991}).
\btitle{The variational form of certain {B}ayes estimators}.
\bjournal{Ann. Statist.}
\bvolume{19}
\bpages{1163--1190}.
\bid{doi={10.1214/aos/1176348244}, issn={0090-5364}, mr={1126320}}
\bptok{imsref}%
\end{barticle}
%
\endbibitem

\bibitem{he2}
%
\begin{bbook}[mr]
\bauthor{\bsnm{Helgason},~\bfnm{Sigurdur}\binits{S.}}
(\byear{1978}).
\btitle{Differential Geometry, {L}ie Groups, and Symmetric Spaces}.
\bseries{Pure and Applied Mathematics}
\bvolume{80}.
\bpublisher{Academic Press}, \baddress{New York}.
\bid{mr={0514561}}
\bptok{imsref}%
\end{bbook}
%
\endbibitem

\bibitem{johnstone08}
%
\begin{barticle}[mr]
\bauthor{\bsnm{Johnstone},~\bfnm{Iain~M.}\binits{I.~M.}}
(\byear{2008}).
\btitle{Multivariate analysis and {J}acobi ensembles: Largest eigenvalue,
{T}racy--{W}idom limits and rates of convergence}.
\bjournal{Ann. Statist.}
\bvolume{36}
\bpages{2638--2716}.
\bid{doi={10.1214/08-AOS605}, issn={0090-5364}, mr={2485010}}
\bptok{imsref}%
\end{barticle}
%
\endbibitem

\bibitem{kkkoz}
%
\begin{bmisc}[auto:STB|2011/12/28|12:52:23]
\bauthor{\bsnm{Kim},~\bfnm{K.~R.}\binits{K.~R.}},
\bauthor{\bsnm{Kim},~\bfnm{P.~T.}\binits{P.~T.}},
\bauthor{\bsnm{Koo},~\bfnm{J.~Y.}\binits{J.~Y.}},
\bauthor{\bsnm{Oh},~\bfnm{J.}\binits{J.}} \AND
\bauthor{\bsnm{Zhu},~\bfnm{H.}\binits{H.}}
(\byear{2011}).
\bhowpublished{An analysis of diffusion tensor image data based on
mixture of
Wisharts. Preprint, Korea Univ.}
\bptok{imsref}%
\end{bmisc}
%
\endbibitem

\bibitem{kr2}
%
\begin{bincollection}[auto:STB|2011/12/28|12:52:23]
\bauthor{\bsnm{Kim},~\bfnm{P.~T.}\binits{P.~T.}} \AND
\bauthor{\bsnm{Richards},~\bfnm{D.~St.~P.}\binits{D.~S.~P.}}
(\byear{2010}).
\btitle{Deconvolution density estimation on spaces of positive definite
symmetric matrices}.
In \bbooktitle{Festschrift for T.~P.~Hettmansperger}
(\beditor{D.~Hunter et al.}, eds.)
\bpages{147--168}.
\bpublisher{World Scientific}, \baddress{Singapore}.
\bptok{imsref}%
\end{bincollection}
%
\endbibitem

\bibitem{Koo}
%
\begin{barticle}[mr]
\bauthor{\bsnm{Koo},~\bfnm{Ja-Yong}\binits{J.-Y.}}
(\byear{1993}).
\btitle{Optimal rates of convergence for nonparametric statistical inverse
problems}.
\bjournal{Ann. Statist.}
\bvolume{21}
\bpages{590--599}.
\bid{doi={10.1214/aos/1176349138}, issn={0090-5364}, mr={1232506}}
\bptok{imsref}%
\end{barticle}
%
\endbibitem

\bibitem{lm}
%
\begin{barticle}[mr]
\bauthor{\bsnm{Letac},~\bfnm{G{\'e}rard}\binits{G.}} \AND
\bauthor{\bsnm{Massam},~\bfnm{H{\'e}l{\`e}ne}\binits{H.}}
(\byear{2007}).
\btitle{Wishart distributions for decomposable graphs}.
\bjournal{Ann. Statist.}
\bvolume{35}
\bpages{1278--1323}.
\bid{doi={10.1214/009053606000001235}, issn={0090-5364}, mr={2341706}}
\bptok{imsref}%
\end{barticle}
%
\endbibitem

\bibitem{linperl}
%
\begin{bincollection}[mr]
\bauthor{\bsnm{Lin},~\bfnm{Shang~P.}\binits{S.~P.}} \AND
\bauthor{\bsnm{Perlman},~\bfnm{Michael~D.}\binits{M.~D.}}
(\byear{1985}).
\btitle{A {M}onte {C}arlo comparison of four estimators of a covariance
matrix}.
In \bbooktitle{Multivariate Analysis {VI} ({P}ittsburgh, {P}A, 1983)}
(\beditor{P. R. Krishnaiah}, ed.)
\bpages{411--429}.
\bpublisher{North-Holland}, \baddress{Amsterdam}.
\bid{mr={0822310}}
\bptok{imsref}%
\end{bincollection}
%
\endbibitem

\bibitem{loh}
%
\begin{barticle}[mr]
\bauthor{\bsnm{Loh},~\bfnm{Wei-Liem}\binits{W.-L.}}
(\byear{1991}).
\btitle{Estimating covariance matrices}.
\bjournal{Ann. Statist.}
\bvolume{19}
\bpages{283--296}.
\bid{doi={10.1214/aos/1176347982}, issn={0090-5364}, mr={1091851}}
\bptok{imsref}%
\end{barticle}
%
\endbibitem

\bibitem{mr}
%
\begin{barticle}[mr]
\bauthor{\bsnm{Mair},~\bfnm{Bernard~A.}\binits{B.~A.}} \AND
\bauthor{\bsnm{Ruymgaart},~\bfnm{Frits~H.}\binits{F.~H.}}
(\byear{1996}).
\btitle{Statistical inverse estimation in {H}ilbert scales}.
\bjournal{SIAM J. Appl. Math.}
\bvolume{56}
\bpages{1424--1444}.
\bid{doi={10.1137/S0036139994264476}, issn={0036-1399}, mr={1409127}}
\bptok{imsref}%
\end{barticle}
%
\endbibitem

\bibitem{Muirhead}
%
\begin{bbook}[mr]
\bauthor{\bsnm{Muirhead},~\bfnm{Robb~J.}\binits{R.~J.}}
(\byear{1982}).
\btitle{Aspects of Multivariate Statistical Theory}.
\bpublisher{Wiley}, \baddress{New York}.
\bid{mr={0652932}}
\bptok{imsref}%
\end{bbook}
%
\endbibitem

\bibitem{pensky}
%
\begin{barticle}[mr]
\bauthor{\bsnm{Pensky},~\bfnm{Marianna}\binits{M.}}
(\byear{1999}).
\btitle{Nonparametric empirical {B}ayes estimation of the matrix
parameter of
the {W}ishart distribution}.
\bjournal{J. Multivariate Anal.}
\bvolume{69}
\bpages{242--260}.
\bid{doi={10.1006/jmva.1998.1803}, issn={0047-259X}, mr={1703374}}
\bptnote{check year}%
\bptok{imsref}%
\end{barticle}
%
\endbibitem

\bibitem{Richards1}
%
\begin{barticle}[mr]
\bauthor{\bsnm{Richards},~\bfnm{Donald St.~P.}\binits{D.~S.~P.}}
(\byear{1985}).
\btitle{Applications of invariant differential operators to multivariate
distribution theory}.
\bjournal{SIAM J. Appl. Math.}
\bvolume{45}
\bpages{280--288}.
\bid{doi={10.1137/0145015}, issn={0036-1399}, mr={0781108}}
\bptok{imsref}%
\end{barticle}
%
\endbibitem

\bibitem{stein-rietz}
%
\begin{bmisc}[auto:STB|2011/12/28|12:52:23]
\bauthor{\bsnm{Stein},~\bfnm{C.}\binits{C.}}
(\byear{1975}).
\bhowpublished{Estimation of a covariance matrix. Rietz lecture, 1975 Annual
Meeting of the IMS, Atlanta, GA}.
\bptok{imsref}%
\end{bmisc}
%
\endbibitem

\bibitem{stein}
%
\begin{bmisc}[mr]
\bauthor{\bsnm{Ste{\u\i}n},~\bfnm{{\v{C}}.}\binits{{\v{C}}.}}
(\byear{1977}).
\bhowpublished{Lectures on the theory of
estimation of many parameters. In \textit{Studies in the Statistical
Theory of Estimation}, \textit{Part I}
(I. A. Ibragimov and M. S. Nikulin, eds.).
\textit{Zap. Nau\v cn. Sem. Leningrad. Otdel. Mat. Inst. Steklov.}
(\textit{LOMI}) \textbf{74} 4--65, 146, 148 (in Russian).
Steklov Mathematical Institute, Moscow.
[English transl.: \textit{J. Soviet Math.} \textbf{34}
(1986) 1373--1403.]}
\bptok{imsref}%
\end{bmisc}
%
\endbibitem

\bibitem{takemura}
%
\begin{barticle}[mr]
\bauthor{\bsnm{Takemura},~\bfnm{Akimichi}\binits{A.}}
(\byear{1984}).
\btitle{An orthogonally invariant minimax estimator of the covariance
matrix of
a multivariate normal population}.
\bjournal{Tsukuba J. Math.}
\bvolume{8}
\bpages{367--376}.
\bid{issn={0387-4982}, mr={0767967}}
\bptok{imsref}%
\end{barticle}
%
\endbibitem

\bibitem{Terras-II}
%
\begin{bbook}[mr]
\bauthor{\bsnm{Terras},~\bfnm{Audrey}\binits{A.}}
(\byear{1988}).
\btitle{Harmonic Analysis on Symmetric Spaces and Applications. {II}}.
\bpublisher{Springer}, \baddress{Berlin}.
\bid{mr={0955271}}
\bptok{imsref}%
\end{bbook}
%
\endbibitem

\bibitem{Zhang}
%
\begin{barticle}[mr]
\bauthor{\bsnm{Zhang},~\bfnm{Cun-Hui}\binits{C.-H.}}
(\byear{1990}).
\btitle{Fourier methods for estimating mixing densities and distributions}.
\bjournal{Ann. Statist.}
\bvolume{18}
\bpages{806--831}.
\bid{doi={10.1214/aos/1176347627}, issn={0090-5364}, mr={1056338}}
\bptok{imsref}%
\end{barticle}
%
\endbibitem

\end{thebibliography}
\end{document}